\documentclass[12pt]{article}
\usepackage{amsfonts}
\usepackage{latexsym, amssymb, amsmath, amscd, amsfonts, epsfig, graphicx}
\usepackage{mathrsfs}
\usepackage{color}
\usepackage{ifpdf}
\newtheorem{theorem}{Theorem}[section]

\newtheorem{example}[theorem]{Example}

 \numberwithin{equation}{section}

\def\qed{\hfill \rule{4pt}{7pt}}

\def\pf{\noindent {\it Proof.} }

\begin{document}

\begin{center}
{\large {\bf Hook Length Formulas for Trees by  Han's Expansion}}

\vskip 6mm

{\small William Y.C. Chen$^1$, Oliver X.Q. Gao$^2$
and Peter L. Guo$^3$\\[%
2mm] Center for Combinatorics, LPMC-TJKLC\\
Nankai University, Tianjin 300071,
P.R. China \\[3mm]
$^1$chen@nankai.edu.cn, $^2$oliver@cfc.nankai.edu.cn, $^3$lguo@cfc.nankai.edu.cn \\[0pt%
] }
\end{center}

\begin{abstract}

Recently Han obtained a general formula for the weight function
corresponding to the expansion  of a generating function in terms of
hook lengths of binary trees. In this paper, we present  formulas
for $k$-ary trees, plane trees, plane forests, labeled trees and
forests. We also find appropriate generating functions which lead to
unifications of the  hook length formulas due to Du and Liu, Han,
Gessel and Seo, and Postnikov.

\end{abstract}
\vskip 3mm


\section{Introduction}

Recently, Han developed an  expansion technique for deriving hook
length formulas for binary trees. He has shown that given any formal
power series $f(x)$ with $f(0)=1$, one can determine the weight
function $\rho(n)$ that leads to a hook length formula for binary
trees. In this paper, we extend Han's technique and obtain the
expansion formulas for $k$-ary trees, plane trees, plane forests,
labeled trees and forests. We find appropriate generating functions
that can be used to derive new hook length formulas, some of which
can be viewed as unifications of the formulas due to Du and Liu
\cite{Du}, Han \cite{Han1,Han2,Han3}, Gessel and Seo \cite{Gel1}.

Let us give a quick review of the background and terminology. For a
tree (or forest) $T$ the hook length of a vertex $u$ of $T$, denoted
by $h_u$, is the number of descendants of $u$ in $T$ under the
assumption that $u$ is counted as a descendant of itself. The hook
length multi-set $\mathcal {H}(T)=\{h_u\colon u\in T\}$ of T is
defined to be the multi-set of hook lengths of the vertices $u$ of
$T$.
\begin{figure}[h,t]
\setlength{\unitlength}{0.5mm}
\begin{center}
\begin{picture}(110,50)
\put(0,15){\circle*{2}}\put(0,30){\circle*{2}}\put(10,45){\circle*{2}}
\put(10,15){\circle*{2}}\put(20,30){\circle*{2}}\put(30,15){\circle*{2}}
\put(0,15){\line(0,1){15}}\put(0,30){\line(2,3){10}}
\put(10,45){\line(2,-3){10}}
\put(10,15){\line(2,3){10}}\put(20,30){\line(2,-3){10}}
\put(10,0){\shortstack{$T$}}\put(50,25){$\mathcal
{H}(T)=\{1,1,1,2,3,6\}$}
\end{picture}
\end{center}
\end{figure}

 Postnikov \cite{Pos} discovered the following remarkable hook length formula
 for binary trees
\begin{equation}{\label{19}}
 \frac{n!}{2^n}\sum_{T}\prod_{h\in
 \mathcal{H}(T)}\left(1+\frac{1}{h}\right)=(n+1)^{n-1},
\end{equation} where the sum ranges over  binary trees with $n$
vertices. Combinatorial proofs of (\ref{19}) have been given by Chen
and Yang \cite{Che3}, and Seo \cite{Seo}. Hook length formulas have
been found for $k$-ary trees, plane forests and forests.  Du and Liu
\cite{Du} have obtained the following formulas
\begin{align}\label{13}
\sum_{T}&\prod_{h\in \mathcal {H}(T)}\left(a+\frac{1}{h}\right)=
\frac{(a+1)}{n!}\prod_{i=1}^{n-1}\big(kan+a+1-i(a-k+1)\big)
\end{align}
and
\begin{align}\label{26}
\sum_{F}\prod_{h\in \mathcal{H}(F)}
\left(a+\frac{1}{h}\right)=\frac{(a+1)}{n!}\prod_{i=1}^{n-1}\big((2n+1)(a+1)-(a+2)i\big),
\end{align}where $T$ (resp., $F$) ranges over  $k$-ary trees (resp., plane forests)
with $n$ vertices. Liu \cite{Liu} has given a hook length formula
for plane forests with a given degree sequence.  Gessel and Seo
\cite{Gel1} independently discovered (\ref{13}) and obtained the
following formula for forests
\begin{align}\label{28}
\sum_{F}\prod_{h\in \mathcal{H}(F)}\left(1+\frac{a}{h}\right)=
(a+1)\prod_{i=1}^{n-1}\big((a+1)n-ai\big),
\end{align} where $F$ ranges over forests with $n$ vertices.

 Han \cite{Han3} has found the following formula for binary trees
\begin{align}\label{40}
n!\sum_{T}&\prod_{h\in \mathcal
{H}(T)}\frac{\big(z+h\big)^{h-1}}{2h(2z+h-1)^{h-2}}=
z\big(n+z\big)^{n-1},
\end{align}where the sum ranges over binary trees with $n$ vertices.
Recall that the polynomials $z(n+z)^{n-1}$ are the classical Abel
polynomials, see Mullin and Rota \cite{Rota}. The above formula was
first proved by Han \cite{Han3} by induction, and then brought into
the framework of his expansion technique \cite{Han1}.

 Han's expansion technique for binary trees
 can be described as follows.    Denote by
$\mathbb{K}[[x]]$ the set of formal power series over some field
$\mathbb{K}$. Define the weight function
$\rho\colon\mathbb{N}^{+}\longrightarrow \mathbb{K}$ to be a mapping
from the set of positive integers to $\mathbb{K}$. Denote by $B(n)$
the set of  binary trees with $n$ vertices. Han \cite{Han1} has
shown that if the following relation holds
\begin{equation} \label{eq-hk}
1+\sum_{n\geq 1}\left(\sum_{T\in B(n)}\prod_{h\in \mathcal
{H}(T)}\rho(h)\right)x^n=f(x),
\end{equation}
then the weight function $\rho$ is given by
\begin{equation}\label{1}
\rho(n)=\frac{[x^n]f(x )}{[x^{n-1}]f(x)^2},
\end{equation}
where $[x^n]f(x)$ denotes the coefficient of $x^n$ in the formal
power series expansion of $f(x)$. The above formula is called the
expansion formula for binary trees. Note that each $T\in B(n)$
$(n\geq 1)$ can be decomposed into a  triple $(T',T'',u)$, where
$T'\in B(m)$ $(0\leq m\leq n-1)$, $T''\in B(n-1-m)$ and  $u$ is the
root of $T$ with hook length $h_u=n$. Suppose that (\ref{eq-hk})
holds. Then we deduce that
\begin{align*}
[x^n]f(x)&=\rho(n)\sum_{(T',T'')}\prod_{h\in
\mathcal {H}(T')}{\rho(h)}\prod_{h\in \mathcal {H}(T'')}{\rho(h)}\\
&=\rho(n)\sum_{m=0}^{n-1}[x^m]f(x) [x^{n-1-m}]f(x),
\end{align*} which implies (\ref{1}).

For example, let $g(x)$ be defined by the functional equation
$g(x)=\textrm{exp}\{xg(x)\}$ and $f(x)=g(2x)$. Applying the Lagrange
inversion formula (see, Stanley \cite[Chapter 5]{Sta2}), we get
\[ [x^n]xg(x)=\frac{1}{n}[x^{n-1}]\textrm{e}^{nx}=\frac{n^{n-1}}{n!}.
\]  Then we have
\[f(x)=\sum_{n\geq 0}(n+1)^{n-1}\frac{(2x)^n}{n!}.\]
By Han's expansion formula for binary trees and the Lagrange
inversion formula for the expansion of $f(x)^2$, the weight function
corresponding to the generating function $f(x)$ is given by
\[\rho(n)=\frac{[x^n]f(x )}{[x^{n-1}]f(x)^2}=1+\frac{1}{n},\]which
leads to Postnikov's formula (\ref{19}).

This paper is organized as follows. In Section 2, we give a
straightforward extension of Han's expansion formula to $k$-ary
trees, and find a generating function whose corresponding hook
length formula is a unification of several known formulas. In
Section 3, we consider the expansion formulas for plane trees and
plane forests. In Section 4, we present expansion formulas and hook
length formulas for labeled trees and forests. We conclude this
paper with the question of finding combinatorial interpretations of
three hook length formulas for plane forests and labeled forests.

\section{$k$-ary trees}

In this section, we begin with a straightforward extension of Han's
expansion formula for binary trees to $k$-ary trees. The formula of
Yang for $k$-ary trees \cite{Yang} corresponds to the expansion of
$e^x$. Moreover, we find a generating function, which is defined by
a functional equation, which enables us to deduce a hook length
formula with one more parameter $z$ compared with the formula
independently due to Du and Liu \cite{Du}, and Gessel and Seo
\cite{Gel1}.

Recall that a \textit{$k$-ary tree} is an ordered rooted unlabeled
tree where each vertex has exactly $k$ subtrees in linear order
where we allow a subtree being empty. When $k=2$, a $k$-ary tree is
called a binary tree. Let $T_k(n)$ denote the set of  $k$-ary trees
with $n$ vertices.

\begin{theorem}
Suppose that we have the following expansion formula for $k$-ary
trees
\begin{align*}
1+\sum_{n\geq 1}\left(\sum_{T\in T_k(n)}\prod_{h\in \mathcal
{H}(T)}\rho(h)\right)x^n=f(x).
\end{align*}
 Then the weight function $\rho$ is given by
\begin{equation}\label{2}
\rho(n)=\frac{[x^n]f(x )}{[x^{n-1}]f(x)^k}.
\end{equation}
\end{theorem}

\begin{example}
Let $f(x)=e^x$. Applying (\ref{2}) we get
$$\rho(n)=\frac{[x^n]e^x}{[x^{n-1}]e^{kx}}=\frac{1}{nk^{n-1}}.$$
Hence
\begin{equation}\label{7}
\sum_{T\in T_k(n)}\prod_{h\in \mathcal
{H}(T)}\frac{1}{hk^{h-1}}=\frac{1}{n!}.
\end{equation}
\end{example}

The above formula (\ref{7}) was derived by Han \cite{Han2} for
$k=2$. Yang \cite{Yang} showed that (\ref{7}) holds for general $k$.
Probabilistic and combinatorial proofs of (\ref{7}) have been given
by Sagan \cite{Sag}, and Chen, Gao and Guo \cite{Che2},
respectively.

Below is another hook length formula of  Han for binary trees.

\begin{theorem}[\mdseries{Han\cite{Han1}, Theorem 6.8}]
For $n\geq 1$,
\begin{equation}\label{37}
\begin{split}
\sum_{T\in B(n)}&\prod_{h\in \mathcal
{H}(T)}\frac{\prod_{i=1}^{h-1}\big(za+z+(2h-i)a+i\big)}
{2h\prod_{i=1}^{h-2}\big(2za+2z+(2h-2-i)a+i\big)}\\&=
\frac{z(a+1)}{n!}\prod_{i=1}^{n-1}\big(za+z+(2n-i)a+i\big).
\end{split}
\end{equation}
\end{theorem}

The generating function $f(x)$ for the above expansion is given by
the following functional equation
\[
g(x)=(a-1)x\big(1+g(x)\big)^{\frac{2a}{a-1}}\] and
\[ f(x)=\big(1+g(x)\big)^{z\frac{2a}{a-1}}.\]

To extend Han's formula to $k$-ary trees, one needs to find the
appropriate extension of the generating function to general $k$.

\begin{theorem} Let $g(x)$ be defined by the functional equation
 \[
g(x)=(a-k)x\big(1+g(x)\big)^{\frac{k(a-1)}{a-k}},\] and $f(x)$ be
given by
\[ f(x)=\big(1+g(x)\big)^{z\frac{a}{a-k}}.\]
Then the weight function $\rho$ corresponding to the hook length
expansion of $f(x)$ for $k$-ary trees is given by
\begin{equation}\label{20}
\rho(n)=\frac{\prod_{i=1}^{n-1}\big(za+k(a-1)n-i(a-k)\big)}
{kn\prod_{i=1}^{n-2}\big(kza+k(a-1)(n-1)-i(a-k)\big)},
\end{equation}
and hence
\begin{equation}\label{18}
\begin{split}
\sum_{T\in T_k(n)}&\prod_{h\in \mathcal
{H}(T)}\frac{\prod_{i=1}^{h-1}\big(za+k(a-1)h-i(a-k)\big)}
{kh\prod_{i=1}^{h-2}\big(kza+k(a-1)(h-1)-i(a-k)\big)}\\&=
\frac{za}{n!}\prod_{i=1}^{n-1}\big(za+k(a-1)n-i(a-k)\big).
\end{split}
\end{equation}
\end{theorem}

\pf By the Lagrange inversion formula we obtain
\begin{align*}
[x^n]f(x)&=\frac{1}{n}[x^{n-1}]z\frac{a}{a-k}
(1+x)^{z\frac{a}{a-k}-1}(a-k)^n(1+x)^{\frac{k(a-1)n}
{a-k}}\\
&=\frac{za(a-k)^{n-1}}{n!}\prod_{i=0}^{n-2}\left(z\frac{a}{a-k}-1+\frac{k(a-1)n}
{a-k}-i\right)\\
&=\frac{za}{n!}\prod_{i=1}^{n-1}\big(za+k(a-1)n-i(a-k)\big).
\end{align*}
Note that $[x^n]f(x)^k$ can be easily derived from $[x^n]f(x)$ by
substituting $z$ with $kz$. Consequently,
\begin{align*}
[x^n]f(x)^k=\frac{kza}{n!}\prod_{i=1}^{n-1}\big(kza+k(a-1)n-i(a-k)\big).
\end{align*}
By Theorem 2.1, we obtain the weight function (\ref{20}). This
yields  (\ref{18}). \qed

The formula (\ref{18}) can be viewed as a unification of several
known formulas. Setting $k=2$ and substituting $a$ with $a+1$ we
obtain Han's formula (\ref{37}). Setting $z=1$ in (\ref{18}) yields
\begin{align*}
\sum_{T\in T_k(n)}&\prod_{h\in \mathcal
{H}(T)}\left(a-1+\frac{1}{h}\right)=
\frac{a}{n!}\prod_{i=1}^{n-1}\big(k(a-1)n+a-i(a-k)\big),
\end{align*}
which is equivalent to the formula (\ref{13}) derived independently
by Du and Liu \cite{Du}, and  Gessel and Seo \cite{Gel1}. Setting
$a=k$ and $z=1$ in (\ref{18}), we obtain
\begin{align*}
\sum_{T\in T_k(n)}&\prod_{h\in \mathcal
{H}(T)}\left(k-1+\frac{1}{h}\right)=
\frac{k^n\big((k-1)n+1\big)^{n-1}}{n!},
\end{align*}
which is an extension of Postnikov's hook formula (\ref{19}) to
$k$-ary trees. Setting $a=k$ in (\ref{18}), we obtain the following
extension of Han's formula (\ref{40}) to $k$-ary trees.

\begin{theorem} For $n\geq 1$,
\begin{align}\label{22}
{n!}\sum_{T\in T_k(n)}&\prod_{h\in \mathcal
{H}(T)}\frac{\big(z+(k-1)h\big)^{h-1}}{kh\big(kz+(k-1)(h-1)\big)^{h-2}}=
z\big(z+(k-1)n\big)^{n-1}.
\end{align}
\end{theorem}

Recall that the polynomials $z\big(z+(k-1)n\big)^{n-1}$ are also
Abel polynomials (see \cite{Rota}). Setting $a\rightarrow \infty$ in
(\ref{18}) we deduce the following formula.

\begin{theorem} For $n\geq 1$,
\begin{align*} \sum_{T\in T_k(n)}&\prod_{h\in \mathcal
{H}(T)}\frac{\prod_{i=1}^{h-1}(kh+z-i)}{kh\prod_{i=1}^{h-2}\big(kh+k(z-1)-i\big)}
=\frac{z}{n!}\prod_{i=1}^{n-1}(kn+z-i).
\end{align*}
\end{theorem}

\section{Plane trees and plane forests}

In this section, we derive hook length  expansion formulas for plane
trees and plane forests. Then we find certain generating functions
to give hook length formulas. Some known formulas can be brought
into the framework of the expansion technique.

A plane tree is a rooted unlabeled tree in which the subtrees,
assumed to be nonempty, of each vertex are arranged in linear order.
A plane forest is a forest of nonempty plane trees which are
linearly ordered. Let $PT(n)$ (resp., $PF(n)$) denote the set of
plane trees (resp., plane forests) with $n$ vertices.

\begin{theorem} \label{ef-pt}
Suppose that  the following expansion formula holds for plane trees
\begin{equation*}
\sum_{n\geq 1}\left(\sum_{T\in PT(n)}\prod_{h\in \mathcal
{H}(T)}\rho(h)\right)x^n=f(x).
\end{equation*}
Then the weight function $\rho$ is given by
\begin{equation}\label{3}
\rho(n)=\frac{[x^n]f(x )}{[x^{n-1}]\frac{1}{1-f(x)}}, \quad n\geq 1.
\end{equation}

\end{theorem}
\pf For a plane tree $T$ with $n$ $(n\geq 2)$ vertices, we obtain a
$j$-tuple $(T_1,T_2,\ldots,T_j)$ $(j\geq 1)$ by deleting the root
$u$ of $T$, where $T_i\in PT(m_i)$ $(m_i>0)$ and $\sum_{i=1}^j
m_i=n-1$. Let $f(n)=[x^n]f(x)$. Then we have
\begin{align*}
f(n)&=\rho(h_u)\sum_{j\geq
1}\sum_{(T_1,T_2,\ldots,T_j)}\prod_{i=1}^j\prod_{h\in
\mathcal {H}(T_i)}{\rho(h)}\\
&=\rho(n)\sum_{j\geq 1}\sum_{m_1+\cdots+m_j=n-1}f(m_1)f(m_2)\cdots
f(m_j)=\rho(n)[x^{n-1}]\frac{1}{1-f(x)}.
\end{align*}
Since $\rho(1)=f(1)$, we arrive at  (\ref{3}). \qed

We proceed to employ (\ref{3}) to derive hook length formulas for
plane trees.

\begin{example}
Let $f(x)$ be defined by
\[ f'(x)=\frac{1}{1-f(x)}\]  with $f(0)=0$.
Solving this differential equation we obtain
 \[ f(x)=\sum_{n\geq
1}(2n-3)!!\frac{x^n}{n!},\]  where $n!!=n(n-2)(n-4)\cdots$ for $n$
odd and $(-1)!!=1$. By (\ref{3}) we see that $ \rho(n)=\frac{1}{n}$.
Hence we are led to the hook length formula
\begin{equation}\label{10}
n!\sum_{T\in PT(n)}\prod_{h\in \mathcal
{H}(T)}{\frac{1}{h}}=(2n-3)!!.
\end{equation}
\end{example}

For a plane tree $T$ with $n$ vertices, it is well-known  that the
number of ways to label the vertices of $T$ with $\{1, 2,\ldots,
n\}$, such that the labeling of each vertex is less than the
labelings of its descendants, is equal to $\frac{n!}{\prod_{h\in
\mathcal {H}(T)}h}$ (see, e.g., \cite{Gel1}). Note that $(2n-3)!!$
equals the number of increasing plane trees on $n$ vertices, see,
e.g., \cite{Gel2}.

\begin{theorem}
For $n\geq 1$, we have
\begin{equation}\label{21}
n!\sum_{T\in PT(n)}\prod_{h\in \mathcal
{H}(T)}{\left(1-\frac{1}{h}\right)^{h-1}}=(n-1)^{n-1}.
\end{equation}
\end{theorem}
\pf Let
\[ f(x)=\sum_{n\geq 1}(n-1)^{n-1}\frac{x^n}{n!}.\] It is known that
(see, Stanley \cite[P. 43]{Sta2}) \[\frac{1}{1-f(x)}=\sum_{n\geq
0}(n+1)^{n-1}\frac{x^n}{n!}.\] Now, by Theorem \ref{ef-pt} we find
\[ \rho(n)=\frac{(n-1)^{n-1}}{n^{n-1}}=\left(1-\frac{1}{n}\right)^{n-1},\]
which gives (\ref{21}). \qed

We next consider hook length formulas for plane forests.

\begin{theorem}
Suppose that  the following expansion formula holds for plane
forests
\begin{equation*}
1+\sum_{n\geq 1}\left(\sum_{T\in PF(n)}\prod_{h\in \mathcal
{H}(T)}\rho(h)\right)x^n=f(x).
\end{equation*}
Then the weight function $\rho$ is given by
\begin{equation}\label{4}
\rho(n)=-\frac{[x^n]f(x)^{-1}}{[x^{n-1}]f(x)}.
\end{equation}
\end{theorem}
\pf To prove (\ref{4}) we notice that
\begin{equation*}
1+\sum_{n\geq 1}\left(\sum_{T\in PF(n)}\prod_{h\in \mathcal
{H}(T)}\rho(h)\right)x^n=\left(1-\sum_{n\geq 1}\Big(\sum_{T\in
PT(n)}\prod_{h\in \mathcal {H}(T)}\rho(h)\Big)x^n\right)^{-1},
\end{equation*}
that is,
\begin{equation*}
\sum_{n\geq 1}\left(\sum_{T\in PT(n)}\prod_{h\in \mathcal
{H}(T)}\rho(h)\right)x^n=1-\left(1+\sum_{n\geq 1}\Big(\sum_{T\in
PF(n)}\prod_{h\in \mathcal {H}(T)}\rho(h)\Big)x^n\right)^{-1}.
\end{equation*}
Utilizing (\ref{3}), it is easy to check (\ref{4}). This completes
the proof. \qed

 Let us consider the above expansion for the exponential function.
\begin{example}
Let $f(x)=e^x$. Then
$$\rho(n)=\frac{(-1)^{n+1}\frac{1}{n!}}{\frac{1}{(n-1)!}}
=\frac{(-1)^{n+1}}{n}.$$
Therefore,
$$1+\sum_{n\geq 1}\sum_{T \in PF(n)}x^n\prod_{h\in
\mathcal {H}(T)}\frac{(-1)^{h+1}}{h}=e^x.$$ Equating the
coefficients of $x^n$ yields
\begin{equation*}
\sum_{T\in PF(n)}\frac{n!}{\prod_{h\in \mathcal
{H}(T)}{(-1)^{h}}{h}}=(-1)^n,
\end{equation*} which can be restated in terms of plane trees
\begin{equation}\label{11}
\sum_{T\in PT(n)}\frac{n!}{\prod_{h\in \mathcal
{H}(T)}{(-1)^{h}}{h}}=-1.
\end{equation}
\end{example}

Note that the formula (\ref{11}) and a combinatorial proof were
given by Yang \cite{Yang}. The following theorem gives an identity
involving Bernoulli numbers.

\begin{theorem}
Let $B_n$ be the $n$$\mathrm{th}$ Bernoulli  number. Then we have
\begin{equation}\label{hl-b}
\sum_{T \in PF(n)}\prod_{h\in \mathcal
{H}(T)}B_h=\frac{(-1)^n}{(n+1)!}.
\end{equation}

\end{theorem}
\pf Let $f(x)=\frac{e^x-1}{x}$. We have
$$\frac{1}{f(x)}=\frac{x}{e^x-1}=\sum_{n\geq 0}B_n\frac{x^n}{n!},$$  By
(\ref{4}) we have $\rho(n)=-B_n$. Hence
\begin{equation}\label{17}
1+\sum_{n\geq 1}\sum_{T \in PF(n)}(-x)^n\prod_{h\in \mathcal
{H}(T)}B_h=\frac{e^x-1}{x}.
\end{equation}
Equating the coefficients of $x^n$ on both sides of (\ref{17}), we
get (\ref{hl-b}). \qed

 The following theorem for plane forests is
a unification of several hook length formulas.

\begin{theorem} Let $g(x)$ be defined by the functional equation
 \[
g(x)=(a+1)x\big(1+g(x)\big)^{\frac{2a}{a+1}},\] and $f(x)$ be given
by
\[ f(x)=\big(1+g(x)\big)^{\frac{za}{a+1}}.\]
Then the weight function $\rho$ corresponding to the hook length
expansion of $f(x)$ for plane forests is given by
\begin{equation}\label{30}
\rho(n)=\frac{\prod_{i=1}^{n-1}\big((2n-z)a-(a+1)i\big)}
{n\prod_{i=1}^{n-2}\big((2n-2+z)a-(a+1)i\big)}.
\end{equation}
Thus we have
\begin{equation}\label{14}
\begin{split}
\sum_{T\in PF(n)}&\prod_{h\in
\mathcal{H}(T)}\frac{\prod_{i=1}^{h-1}\big((2h-z)a-(a+1)i\big)}
{h\prod_{i=1}^{h-2}\big((2h-2+z)a-(a+1)i\big)}
\\
&=\frac{za}{n!}\prod_{i=1}^{n-1}\big((2n+z)a-(a+1)i\big).
\end{split}
\end{equation}
\end{theorem}

\pf By the  Lagrange inversion formula we get
\begin{align*}
[x^n]f(x)&=\frac{1}{n}[x^{n-1}]\frac{za}{a+1}(1+x)^
{\frac{za}{(a+1)}-1}(a+1)^n(1+x)^{\frac{2an}{a+1}}\\
&=\frac{za}{n!}\prod_{i=1}^{n-1}\big((2n+z)a-(a+1)i\big).
\end{align*}
Substituting $z$ with $-z$ yields
\begin{equation*}
[x^n]f(x)^{-1}
=\frac{-za}{n!}\prod_{i=1}^{n-1}\big((2n-z)a-(a+1)i\big).
\end{equation*}
By Theorem 3.4,  it is easy to verify (\ref{30}) and (\ref{14}).
\qed

We now consider some special cases of the above formula (\ref{14}).
Taking $z=1$ in (\ref{14}) we get the formula (\ref{26}) derived by
Du and Liu \cite{Du}
\begin{equation*}
\sum_{T\in PF(n)}\prod_{h\in \mathcal{H}(T)}
\left(a-1+\frac{1}{h}\right)=\frac{a}{n!}\prod_{i=1}^{n-1}\big((2n+1)a-(a+1)i\big).
\end{equation*}

Taking $a=1$ in (\ref{14}) we have the following  identity.
\begin{theorem}
For $n\geq 1$,
\begin{align}\label{32}
\sum_{T\in PF(n)}&\prod_{h\in
\mathcal{H}(T)}\frac{\prod_{i=1}^{h-1}(2h-z-2i)}
{h\prod_{i=2}^{h-1}(2h+z-2i)}
=\frac{z}{n!}\prod_{i=1}^{n-1}(2n+z-2i).
\end{align}
\end{theorem}

Taking $a=-1$ in (\ref{14}) we deduce the following formula.

\begin{theorem}
For $n\geq 1$,
\begin{equation}\label{23}
\sum_{T\in PF(n)}\prod_{h\in
\mathcal{H}(T)}\frac{(2h-z)^{h-1}}{h(2h-2+z)^{h-2}}=\frac{z}{n!}(2n+z)^{n-1}.
\end{equation}
\end{theorem}

Setting $z=2$ in (\ref{23}) we get the a formula  equivalent to
(\ref{21}).

\begin{theorem}
For $n\geq 1$,
\begin{align}\label{12}
n!\sum_{T\in PF(n)}&\prod_{h\in
\mathcal{H}(T)}\left(1-\frac{1}{h}\right)^{h-1}=(n+1)^{n-1}.
\end{align}
\end{theorem}

Taking $a\rightarrow \infty$ in (\ref{14}) we are led to the
following identity.

\begin{theorem}
For $n\geq 1$,
\begin{equation*}
\sum_{T\in PF(n)}\prod_{h\in
\mathcal{H}(T)}\frac{(2h-z-1)_{h-1}}{h(2h+z-3)_{h-2}}=\frac{z}{n!}(2n+z-1)_{n-1},
\end{equation*}
where $(x)_n=x(x-1)\cdots(x-n+1)$ stands for the falling factorial.
\end{theorem}

\section{Labeled trees and forests}

In this section, we present the expansion formula for labeled trees
and forests. Numerous new hook length formulas are derived. In
particular, we derive a unified hook length formula which includs
the formula (\ref{28}) derived by Gessel and Seo \cite{Gel1} as its
special case. Let $T(n)$ (resp., $F(n)$) denote the set of labeled
trees (resp., forests) with $n$ vertices.

\begin{theorem}
Suppose that the following expansion formula holds for labeled trees
\begin{equation*}
\sum_{n\geq 1}\left(\sum_{T\in T(n)}\prod_{h\in \mathcal
{H}(T)}\rho(h)\right)\frac{x^n}{n!}=f(x).
\end{equation*}
Then the weight function $\rho$ is given by
\begin{equation}\label{5}
\rho(n)=\frac{[x^n]f(x )}{[x^{n-1}]e^{f(x)}}.
\end{equation}
\end{theorem}

\pf For a labeled tree $T$ with $n$ $(n\geq 2)$ vertices, let $u$ be
the root of $T$.  Let $\{T_1,T_2,\ldots,T_j\}$ $(j\geq 1)$ be the
set of subtrees of $u$, and let $B_i$ be the underlying set of
$T_i$. Keep in mind that $\{B_1,\ldots,B_j\}$ is a  partition of
$[n]-\{u\}$. Put  $f(n)=[x^n]f(x)$. Since there are $n$ choices for
$u$, we find that
\begin{equation*}
f(n)=\rho(n)n\sum_{j\geq 1}\sum_{\{B_1,\ldots,B_j\}}f(\# B_1)\cdots
f(\# B_j),
\end{equation*} where $\{B_1,\ldots,B_j\}$ ranges over all set partitions
of $[n-1]$ with $j$ blocks and $\# B_i$ denotes the cardinality of
$B_i$. From the exponential formula  (see, e.g., Stanley
\cite{Sta2}) it follows that
$$\sum_{j\geq 1}\sum_{\{B_1,\ldots,B_j\}}f(\#
B_1)\cdots f(\# B_j)=\left[\frac{x^{n-1}}{(n-1)!}\right]e^{f(x)}.$$
Since $\rho(1)=f(1)$, we conclude that  (\ref{5}) is true for $n\geq
1$. \qed

\begin{theorem}
Let $B_n$ be the $n$-$\mathrm{th}$ Bell number; that is, the number
of partitions of $[n]$. Then we have$$\sum_{T\in T(n)}\prod_{h\in
\mathcal{H}(T)}\frac{1}{hB_{h-1}}=1.$$
\end{theorem}
\pf Let $f(x)=e^x-1$. Then \[ e^{f(x)}=\sum_{n\geq 0}B_n
\frac{x^n}{n!}.\] It follows from (\ref{5}) that \[
\rho(n)=\frac{1}{nB_{n-1}},\] as desired. \qed

\begin{theorem}
For $n\geq 1$, we have
\begin{align}\label{16}
\sum_{T\in T(n)}\prod_{h\in \mathcal{H}(T)}{\frac{1}{h}}=(n-1)!.
\end{align}
\end{theorem}

\pf Let $$f(x)=\ln\frac{1}{1-x}.$$ Clearly,
$e^{f(x)}=\frac{1}{1-x}$. By (\ref{5}) we see that
$\rho(n)=\frac{1}{n}$. This completes the proof. \qed

The next theorem gives the hook length expansion formula for
forests.

\begin{theorem}
Suppose that  the following expansion formula holds
\begin{equation*}
1+\sum_{n\geq 1}\left(\sum_{T\in F(n)}\prod_{h\in \mathcal
{H}(T)}\rho(h)\right)\frac{x^n}{n!}=f(x).
\end{equation*}
 Then the weight function $\rho$ is given by
\begin{equation}\label{6}
\rho(n)=\frac{[x^n]\ln{f(x )}}{[x^{n-1}]f(x)}.
\end{equation}
\end{theorem}

\pf From the the exponential formula (see, Stanley \cite[Chpter
5]{Sta2}), we see that
\begin{equation*}
f(x)=1+\sum_{n\geq 1}\left(\sum_{T\in F(n)}\prod_{h\in \mathcal
{H}(T)}\rho(h)\right)\frac{x^n}{n!}=\mathrm{exp}\left\{\sum_{n\geq
1}\Big(\sum_{T\in T(n)}\prod_{h\in \mathcal
{H}(T)}\rho(h)\Big)\frac{x^n}{n!}\right\}.
\end{equation*}
Taking logarithm on both sides yields
\begin{equation}\label{41}
\ln f(x) =\sum_{n\geq 1}\Big(\sum_{T\in T(n)}\prod_{h\in \mathcal
{H}(T)}\rho(h)\Big)\frac{x^n}{n!}.
\end{equation}
Combining the above relation and the expansion formula (\ref{5}) for
labeled rooted trees,  we arrive at (\ref{6}). \qed

\begin{theorem}
For $n\geq 1$,
$$\frac{1}{n!}\sum_{T\in F(n)}\prod_{h\in
\mathcal{H}(T)}\frac{2(2h-2)!!}{h(2h-3)!!}={2n \choose n}.$$
\end{theorem}

\pf Let \[ f(x)=\frac{1}{\sqrt{1-4x}}=\sum_{n\geq 0}{2n \choose
n}x^n.\] Applying  (\ref{6}) we obtain \[
\rho(n)=\frac{2(2n-2)!!}{n(2n-3)!!},\] as desired. \qed

The following formula is a unification of several hook length
formulas  for forests.

\begin{theorem} Let $g(x)$ be defined by the functional equation
 \[
g(x)=x(a-1)\big(1+g(x)\big)^{\frac{a}{a-1}},\] and $f(x)$ be given
by
\[f(x)=\big(1+g(x)\big)^{\frac{za}{a-1}}.\]
Then the weight function $\rho$ corresponding to the hook length
expansion of $f(x)$ for forests is given by
\begin{equation}\label{31}
\rho(n)=\frac{\prod_{i=1}^{n-1}\big(an-(a-1)i\big)}
{n\prod_{i=1}^{n-2}\big(a(n-1+z)-(a-1)i\big)}.
\end{equation}
So we have
\begin{equation}\label{25}
\begin{split}
\sum_{T\in F(n)}&\prod_{h\in
\mathcal{H}(T)}\frac{\prod_{i=1}^{h-1}\big(ah-(a-1)i\big)}
{h\prod_{i=1}^{h-2}\big(a(h-1+z)-(a-1)i\big)}\\
&=za\prod_{i=1}^{n-1}\big(a(n+z)-(a-1)i\big).
\end{split}
\end{equation}

\end{theorem}
\pf By the Lagrange inversion formula we get
\begin{align*}
[x^n]f(x)&=\frac{1}{n}[x^{n-1}]\frac{za}{a-1}(1+x)^{\frac{za}{a-1}-1}
(a-1)^n(1+x)^{\frac{an}{a-1}}\\
&=\frac{za}{n!}\prod_{i=1}^{n-1}\big(a(n+z)-(a-1)i\big)
\end{align*}and
\begin{align*}
[x^n]\ln f(x)&=\frac{1}{n}[x^{n-1}]\frac{za}{a-1}(1+x)^{-1}
(a-1)^n(1+x)^{\frac{an}{a-1}}\\
&=\frac{za}{n!}\prod_{i=1}^{n-1}\big(an-(a-1)i\big).
\end{align*}
By Theorem 4.4, it is easy to verify (\ref{31}), and so we obtain
(\ref{25}). \qed

Setting $z=1$ in (\ref{25}) we obtain the following formula
equivalent to (\ref{28}) derived by Gessel and Seo \cite{Gel1}:
\begin{equation*}
\sum_{T\in F(n)}\prod_{h\in
\mathcal{H}(T)}\left(1+\frac{a-1}{h}\right)=
a\prod_{i=1}^{n-1}\big(an-(a-1)i\big).
\end{equation*}

Setting $a=1$ in (\ref{25}) we have the following formula.
\begin{theorem}
For $n\geq 1$,
\begin{equation}\label{34}
\sum_{T\in F(n)}\prod_{h\in
\mathcal{H}(T)}\Big(\frac{h}{h-1+z}\Big)^{h-2}=z(n+z)^{n-1}.
\end{equation}
\end{theorem}

Taking $a\rightarrow \infty$ in (\ref{25}) leads to the following
identity.

\begin{theorem}
For $n\geq 1$,
\begin{equation}\label{35}
\sum_{T\in F(n)}\prod_{h\in
\mathcal{H}(T)}\frac{(h-1)!}{h(h-2+z)_{h-2}}=z(n+z-1)_{n-1}.
\end{equation}
\end{theorem}

Setting $z=1$ in (\ref{35}) yields a formula which is  equivalent to
(\ref{16}).

\begin{theorem}
For $n\geq 1$,
\begin{equation}\label{38}
\sum_{T\in F(n)}\prod_{h\in \mathcal{H}(T)}\frac{1}{h}=n!.
\end{equation}
\end{theorem}

Setting $z=2$ in (\ref{35}) we get the following  formula.

\begin{theorem}
For $n\geq 1$,
\begin{equation}\label{36}
\sum_{T\in F(n)}\prod_{h\in
\mathcal{H}(T)}\frac{1}{h^2}=\frac{(n+1)!}{2^{n}}.
\end{equation}
\end{theorem}

It would be interesting to find combinatorial proofs of the formulas
(\ref{12}), (\ref{38}) and (\ref{36}). Notice that the right hand
side of (\ref{12}) is  equal to the number of forests with $n$
vertices and for a forest $F$,
\[ n!\over \prod_{h\in
\mathcal{H}(F)} h
\] equals the number of increasing labelings of $F$, see Gessel and Seo \cite{Gel1}.

\vspace{0.5cm}
 \noindent{\bf Acknowledgments.}  We would like to thank Guo-Niu Han for
 valuable suggestions.  This work was supported by  the 973
Project, the PCSIRT Project of the Ministry of Education, the
Ministry of Science and Technology, and the National Science
Foundation of China.


\end{document}